\documentclass{amsart}%
\pdfoutput=1
\usepackage{amsfonts}
\usepackage{amsmath}
\usepackage{amssymb}
\usepackage{graphicx}%
\providecommand{\U}[1]{\protect\rule{.1in}{.1in}}
\pdfoutput=1
\newtheorem{theorem}{Theorem}
\theoremstyle{plain}

\numberwithin{equation}{section}
\begin{document}
\title[Odd-Graceful Labelings]{Odd-Graceful Labelings of Trees of Diameter 5}
\author{Christian Barrientos}
\address{Department of Mathematics\\
Clayton State University\\
Morrow, GA 30260, USA}
\email{chr\_barrientos@yahoo.com}
\urladdr{http://cims.clayton.edu/cbarrien}
\thanks{This paper is in final form and no version of it will be submitted for
publication elsewhere.}

\begin{abstract}
A difference vertex labeling of a graph $G$ is an assignment $f$ of labels to
the vertices of $G$ that induces for each edge $xy$ the weight $\left\vert
f(x)-f(y)\right\vert .$ A difference vertex labeling $f$ of a graph $G$ of
size $n$ is odd-graceful if $f$ is an injection from $V(G)$ to
$\{0,1,...,2n-1\}$ such that the induced weights are $\{1,3,...,2n-1\}.$ We
show here that any forest whose components are caterpillars is odd-graceful.
We also show that every tree of diameter up to five is odd-graceful.

\end{abstract}
\maketitle

\section{Introduction}

Let $G$ be a graph of order $m$ and size $n,$ a \textit{difference vertex
labeling} of $G$ is an assignment $f$ of labels to the vertices of $G$ that
induces for each edge $xy$ a label or \textit{weight} given by the absolute
value of the difference of its vertex labels. Graceful labelings are a
well-known type of difference vertex labeling; a function $f$ is a
\textit{graceful labeling} of a graph $G$ of size $n$ if $f$ is an injection
from $V(G)$ to the set $\{0,1,...,n\}$ such that, when each edge $xy$ of $G$
has assigned the weight $\left\vert f(x)-f(y)\right\vert ,$ the resulting
weights are distinct; in other words, the set of weights is $\{1,2,...,n\}$. A
graph that admits a graceful labeling is said to be \textit{graceful}.

When a graceful labeling $f$ of a graph $G$ has the property that there exists
an integer $\lambda$ such that for each edge $xy$ of $G$ either $f(x)\leq
\lambda<f(y)$ or $f(y)\leq\lambda<f(x),$ $f$ is named an $\alpha$-labeling and
$G$ is said to be an $\alpha$-graph. From the definition it is possible to
deduce that an $\alpha$-graph is necessarily bipartite and that the number
$\lambda$ (called the \textit{boundary value} of $f$) is the smaller of the
two vertex labels that yield the edge with weight 1. Some examples of $\alpha
$-graphs are the cycle $C_{n}$ when $n\equiv0(\operatorname{mod}4),$ the
complete bipartite graph $K_{m,n},$ and caterpillars (i.e., any tree with the
property that the removal of its end vertices leaves a path).

A little less restrictive than $\alpha$-labelings are the odd-graceful
labelings introduced by Gnanajothi in 1991 \cite{4}. A graph $G$ of size $n$
is \textit{odd-graceful} if there is an injection $f:V(G)\rightarrow
\{0,1,2,...,2n-1\}$ such that the set of induced weights is
$\{1,3,...,2n-1\}.$ In this case, $f$ is said to be an \textit{odd-graceful
labeling} of $G.$ One of the applications of these labelings is that trees of
size $n,$ with a suitable odd-graceful labeling, can be used to generate
cyclic decompositions of the complete bipartite graph $K_{n,n}.$ In Figure 1
we show an odd-graceful tree of size 6 together with its embedding in the
circular arrangement used to produce the cyclic decomposition of $K_{6,6}.$
Once the labeled tree has been embedded, succesives $30^{\circ}$
(counterclockwise) rotations produce the desired cyclic decomposition of
$K_{6,6}.$%

\[%
{\parbox[b]{2.6426in}{\begin{center}
\includegraphics[
natheight=1.398400in,
natwidth=2.601300in,
height=1.4338in,
width=2.6426in
]%
{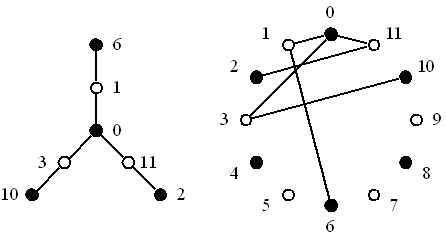}%
\\
Fig. 1. Cyclic decomposition of $K_6,6$
\end{center}}}%
\]

Gnanajothi \cite{4} proved that the class of odd-graceful graphs lies between
the class of $\alpha$-graphs and the class of bipartite graphs; she proved
that every $\alpha$-graph is also odd-graceful. The reverse case does not
work, for example the odd-graceful tree shown in Figure 1 is the smallest tree
without an $\alpha$-labeling. Since many families of $\alpha$-graphs are
known, the most attractive examples of odd-graceful graphs are those without
an $\alpha$-labeling or where an $\alpha$-labeling is unknown; for instance,
Gnanajothi [4] proved that the following are odd-graceful graphs: $C_{n}$ when
$n\equiv2(\operatorname{mod}4),$ the disjoint union of $C_{4},$ the prism
$C_{n}\times K_{2}$ if and only if $n$ is even, and trees of diameter 4 among
others. Eldergil \cite{2} proved that the one-point union of any number of
copies of $C_{6}$ is odd-graceful. Seoud, Diab, and Elsakhawi \cite{5} showed
that a connected $n$-partite graph is odd-graceful if and only if $n=2$ and
that the join of any two connected graphs is not odd-graceful.

A detailed account of results in the subject of graph labelings can be found
in Gallian' survey \cite{3}.

Gnanajothi \cite{4} conjectured that all trees are odd-graceful and verified
this conjecture for all trees with order up to 10. The author has extended
this up to trees with order up to 12\footnote{Odd-graceful labelings of trees
of order 11 and 12 can be found at http://cims.clayton.edu/cbarrien/research}.
In this paper we prove that all trees of diameter 5 are odd-graceful and that
any forest whose components are caterpillars is odd-graceful.

\section{Odd-Graceful Forests}

In this section we study $forests$ that accept odd-graceful labelings. Recall
that a forest with more than one component cannot be graceful bucause it has
"too many edges". First we prove that any graph that admits an $\alpha
$-labeling also admits an odd-graceful labeling by transforming conviniently
its $\alpha$-labeling.

\begin{theorem}
Any $\alpha$-graph is odd-graceful.
\end{theorem}

\begin{proof}
Let $G$ be an $\alpha$-graph of size $n,$ as consequence $G$ is bipartite with
partition $\{A,B\}.$ Suppose that $f$ is an $\alpha$-labeling of $G$ such that
$\max\{f(x):x\in A\}<\min\{f(x):x\in B\}.$ Let $g$ be a labeling of the
vertices of $G$ defined by
\[
g(x)=\left\{
\begin{array}
[c]{ll}%
2f(x), & x\in A\\
2f(x)-1, & x\in B.
\end{array}
\right.
\]
Thus, the labels assigned by $g$ are in the set $\{0,1,...,2n-1\},$
furthermore, the weight of the edge $xy$ of $G$ induced by the labeling $f,$
where $x\in A$ and $y\in B,$ is $w=f(y)-f(x),$ so its weight under the
labeling $g$ is $g(y)-g(x)=2f(y)-1-2f(x)=2(f(y)-f(x))-1=2w-1.$ Since $1\leq
w\leq n,$ we have that the weights induced by $g$ are $\{1,3,...,2n-1\}.$
Therefore, $g$ is an odd-graceful labeling of $G.$
\end{proof}

In Figure 2 we show an example of an $\alpha$-labeling of a caterpillar,
followed for the corresponding odd-graceful labeling. We use this labeling in
the next theorem.%

\[%
{\parbox[b]{3.8233in}{\begin{center}
\includegraphics[
natheight=0.773100in,
natwidth=3.773200in,
height=0.8043in,
width=3.8233in
]%
{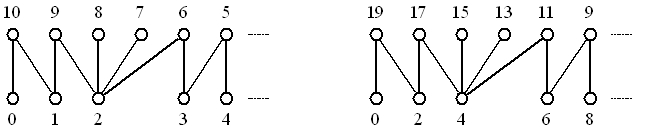}%
\\
Fig. 2. Odd-graceful labeling of a caterpillar
\end{center}}}%
\]

\begin{theorem}
Any forest which components are caterpillars is odd-graceful.
\end{theorem}

\begin{proof}
Let $F_{i}$ be a caterpillar of size $n_{i}\geq1,$ for $1\leq i\leq k.$ Let
$u_{i},v_{i}\in V(F_{i})$ such that $d(u_{i},v_{i})=diam(F_{i});$ so
identifying $v_{i}$ with $u_{i+1},$ for each $1\leq i\leq k-1,$ we have a
caterpillar $F$ of size $\sum\limits_{i=1}^{k}n_{i}=n.$ Now we proceed to find
both, the $\alpha$-labeling of $F$ and its corresponding odd-graceful
labeling, using the scheme shown in Figure 2. Once the odd-graceful labeling
has been obtained, we disengage each caterpillar $F_{i}$ from $F,$ keeping
their labels; in this form, the weights induced are $\{1,3,...,2n-1\}.$ To
eliminate the overlapping of labels we subtract 1 from each vertex label of
$F_{i}$ when $i$ is even, in this way the weights remain the same and the
labels assigned on $u_{i+1}and$ $v_{i}$ differ by one unit. Therefore, the
labeling of the forest $\bigcup\limits_{i=1}^{k}F_{i}$ is odd-graceful.
\end{proof}

In Figure 3 we show an example of this construction using the odd-graceful
labeling obtained in Figure 2.%

\[%
{\parbox[b]{4.0733in}{\begin{center}
\includegraphics[
natheight=1.718400in,
natwidth=4.023100in,
height=1.7564in,
width=4.0733in
]%
{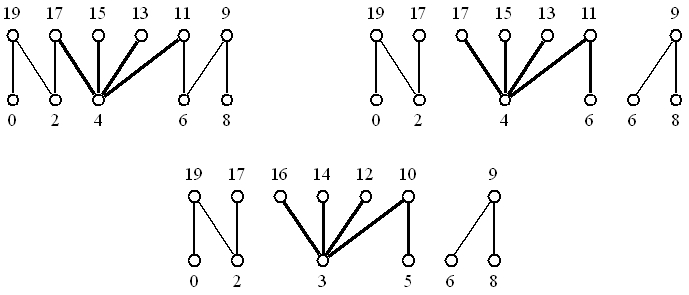}%
\\
Fig. 3. Odd-graceful labeling of a caterpillar
\end{center}}}%
\]

The procedure used in this proof can be extended to the disjoint union of
graphs with $\alpha$-labelings. In fact, suppose that the concatenation of
blocks $B_{1},B_{2},...,B_{k}$ results in a graph $G$ whose block-cutpoint
graph is a path. In \cite{1} we proved that if each $B_{i}$ is an $\alpha
$-graph, so it is $G.$ Transforming this $\alpha$-labeling into an
odd-graceful labeling and disconnecting $G$ into blocks, the disjoint union of
these blocks is odd-graceful.

\begin{theorem}
The disjoint union of blocks that accept $\alpha$-labelings is odd-graceful.
\end{theorem}

\section{Odd-Graceful Trees of Diameter Five}

Every tree of diameter at most 3 is a caterpillar, therefore it is
odd-graceful. Gnanajothi [4] proved that every rooted tree of height 2 (that
is, diameter 4) is odd-graceful. In the next theorem we represent trees of
diameter 5 as rooted trees of height 3 and prove that they are odd-graceful.

Let $T$ be a tree of diameter 5; $T$ can be represented as a rooted tree of
height 3 by using any of its two central vertices as the root vertex. Note
that only one of the vertices in level 1 has descendants in level 3; this
vertex will be located in the right extreme of level 1. Now, within each
level, the vertices are placed from left to right in such a way that their
degrees are increasing. In the proof of the next theorem we use this type of
representation of $T,$ that is, assuming that $v$ (one of the two central
vertices) is the root.

\begin{theorem}
All trees of diameter five are odd-graceful.
\end{theorem}

\begin{proof}
Let $T$ be a tree of diameter 5 and size $n.$ Suppose that $T$ has been drawn
according to the previous description. Let $v_{i,j}$ denote the $i$th vertex
of level $j,$ for $j=1,2,3,$ this vertex is placed at the right of
$v_{i+1,j}.$ Consider the labeling $f$ of the vertices within each level given
by recurrence as follows: $f(v)=0,$ $f(v_{1,1})=2n-2\deg(v)+1,$ $f(v_{1,2}%
)=2,$ $f(v_{1,3})=3,$ and $f(v_{i,j})=f(v_{i-1,j})+d(v_{i,j},v_{i-1,j})$ where
$i\geq2$ and $1\leq j\leq3.$

We claim that $f$ is an odd-graceful labeling of $T.$ In fact, let us see that
there is no overlapping of labels. On level $0$ the label used is $0$ and on
level $2$ all labels are even being $2$ the smallest label used here. On
levels $1$ and $3$ the labels used are odd; on level $1$ the labels used are
$2n-1,2n-3,...,2n-2\deg(v)+1,$ while on level $3$ the labels used are
$3,5,...,2\deg(v_{1,2})-1.$ Now we need to prove that $2n-2\deg(v)+1>2\deg
(v_{1,2})-1;$ since $T$ is a tree of diameter $5,$ at least two vertices on
level $1$ has descendants, so $n+1>\deg(v)+\deg(v_{1,2}),$ which implies the
desired inequality.

As a consequence of the fact that labels used in consecutive levels have
different parity, each weight obtained is an odd number not exceeding $2n-1.$
Suppose that $v_{i+1,j}$ and $v_{i,j}$ have the same father $x,$ by definition
of $f,$ the edges $xv_{i+1,j}$ and $xv_{i,j}$ have consecutive weights. If
$v_{i+1,j}$ and $v_{i,j}$ have different father, $x$ and $y$ respectively,
then $\left\vert f(y)-f(v_{i,j})\right\vert =\left\vert (f(x)+2)-(f(v_{i+1,j}%
)+4)\right\vert =\left\vert f(x)-f(v_{i+1,j})-2\right\vert .$ Thus, the
weights are $2n-2\deg(v)-1,...,2\deg(v)+1.$ On level 2, the weights are
$2n-2\deg(v)-3,...,2\deg(v_{1,2})-1,$ and on level 3 the weights are
$2\deg(v_{1,2})-3,...,1.$

Therefore, $f$ is an odd-graceful labeling of $T.$
\end{proof}

In Figure 4 we present a scheme of this labeling for a tree of size 13.%

\[%
{\parbox[b]{2.0271in}{\begin{center}
\includegraphics[
natheight=1.523800in,
natwidth=1.988200in,
height=1.5601in,
width=2.0271in
]%
{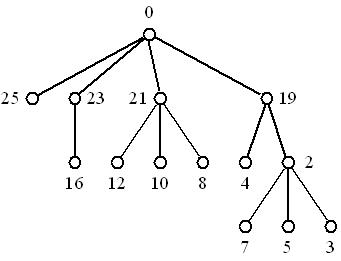}%
\\
Fig. 4. Odd-graceful tree of diameter 5
\end{center}}}%
\]

Similar arguments can be used to find odd-graceful labelings of trees of
diameter 6; however we do not have a general labeling scheme for this case. So
it is an open problem determining whether trees of diameter 6 are
odd-graceful. In Figure 5, we give an example of an odd-graceful labeling for
a tree of size 17 and diameter 6.%

\[%
{\parbox[b]{2.7069in}{\begin{center}
\includegraphics[
natheight=1.523800in,
natwidth=2.664500in,
height=1.5601in,
width=2.7069in
]%
{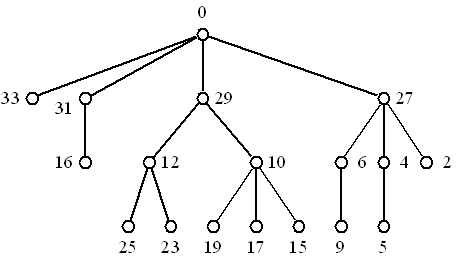}%
\\
Fig. 5. Odd-graceful tree of diameter 6
\end{center}}}%
\]

To conclude this section, we show in Figure 6 an odd-graceful labeling for a
special type of tree of diameter 6, namely the star $S(n,3)$ with $n$ spokes
of length $3.$%

\[%
{\parbox[b]{2.6429in}{\begin{center}
\includegraphics[
natheight=1.460700in,
natwidth=2.601300in,
height=1.4961in,
width=2.6429in
]%
{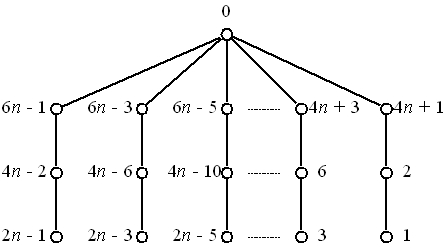}%
\\
Fig. 6. Odd-graceful labeling of the star $S(n,3)$
\end{center}}}%
\]


\begin{thebibliography}{9}                                                                                                %


\bibitem {1}C. Barrientos, Graceful labelings of chain and corona graphs,
\textit{Bull. Inst. Combin. Appl.}, \textbf{34}(2002) 17-26

\bibitem {2}P. Eldergill, Decomposition of the Complete Graph with an Even
Number of Vertices. M. Sc. Thesis, McMaster University, 1997

\bibitem {3}J.A. Gallian, A dynamic survey of graph labeling,
\textit{Electron. J. Combinatorics} (2008), \#DS6

\bibitem {4}R.B. Gnanajothi, Topics in Graph Theory. Ph. D. Thesis. Madurai
Kamaraj University, 1991

\bibitem {5}M.A. Seoud, A.E.I. Abdel Maqsoud, and E.A. Elsahawi, On
strongly-\textit{C} harmoniuos, relatively prime, odd graceful and cordial
graphs, \textit{Proc. Math. Phys. Soc. Egypt}, \textbf{73}(1998) 33-55
\end{thebibliography}
\end{document}